\documentclass[pdflatex,sn-basic]{sn-jnl} 

\usepackage{graphicx}%
\usepackage{multirow}%
\usepackage{amsmath,amssymb,amsfonts}%
\usepackage{amsthm}%
\usepackage{mathrsfs}%
\usepackage[title]{appendix}%
\usepackage{xcolor}%
\usepackage{textcomp}%
\usepackage{manyfoot}%
\usepackage{booktabs}%
\usepackage{algorithm}%
\usepackage{algorithmicx}%
\usepackage{algpseudocode}
\usepackage{listings}%
\usepackage{lipsum}
\usepackage{epstopdf}
\begin{document}

\title[Article Title]{Performance Comparison of Design Optimization and Deep Learning-based Inverse Design}


\author[1]{\fnm{Minyoung} \sur{Jwa}}\email{minyoung306@sookmyung.ac.kr}
\equalcont{These authors contributed equally to this work.}

\author[2]{\fnm{Jihoon} \sur{Kim}}\email{jihoon.kim@kaist.ac.kr}
\equalcont{These authors contributed equally to this work.}

\author[3]{\fnm{Seungyeon} \sur{Shin}}\email{seungyeon@kaist.ac.kr}

\author[1]{\fnm{Ah-hyeon} \sur{Jin}}\email{lisa7399@naver.com}

\author[3]{\fnm{Dongju} \sur{Shin}}\email{dongju1118@kaist.ac.kr}

\author*[3,4]{\fnm{Namwoo} \sur{Kang}}\email{nwkang@kaist.ac.kr}

\affil[1]{\orgdiv{Department of Mechanical Systems Engineering}, \orgname{Sookmyung Women's University}, \orgaddress{\street{100, Cheongpa-ro 47-gil}, \city{Seoul}, \postcode{04310}, \country{Korea}}}

\affil[2]{\orgdiv{Department of Mechanical Engineering}, \orgname{Korea Advanced Institute of Science and Technology}, \orgaddress{\street{291, Daehak-ro}, \city{Daejeon}, \postcode{34141}, \country{Korea}}}

\affil*[3]{\orgdiv{CCS Graduate School of Mobility}, \orgname{Korea Advanced Institute of Science and Technology}, \orgaddress{\street{193, Munji-ro}, \city{Daejeon}, \postcode{34051}, \country{Korea}}}

\affil[4]{\orgname{Narnia Labs}, \orgaddress{\street{193, Munji-ro}, \city{Daejeon}, \postcode{34051}, \country{Korea}}}

\abstract{
    Surrogate model-based optimization has been increasingly used in the field of engineering design.
    It involves creating a surrogate model with objective functions or constraints based on the data obtained from simulations or real-world experiments, and then finding the optimal solution from the model using numerical optimization methods.
    Recent advancements in deep learning-based inverse design methods have made it possible to generate real-time optimal solutions for engineering design problems, eliminating the requirement for iterative optimization processes.
    Nevertheless, no comprehensive study has yet closely examined the specific advantages and disadvantages of this novel approach compared to the traditional design optimization method.
    The objective of this paper is to compare the performance of traditional design optimization methods with deep learning-based inverse design methods by employing benchmark problems across various scenarios.
    Based on the findings of this study, we provide guidelines that can be taken into account for the future utilization of deep learning-based inverse design. It is anticipated that these guidelines will enhance the practical applicability of this approach to real engineering design problems.
    
}


\maketitle

\section{Introduction}\label{sec:intro}

Design optimization of most engineering systems rely on complex and time-consuming physics-based simulations.
Surrogate-based optimization has a potential to be a more efficient approach that replaces such simulations with a surrogate model that gives faster convergence with less computational costs \citep{queipo2005surrogate}, as studied in the field of computational fluid dynamics and computational structrual dynamics \citep{han2012surrogate}.
Recently, deep learning-based (deep) surrogate modeling has become an active research area in the field of design optimization \citep{tang2020subsurface} due to their ability in generating good estimations with less time costs.

However, design optimization relying on traditional optimization algorithms still has inherent limitations in terms of both time and computational costs in order to be applied in practical industrial fields. 
To address this issue, many deep inverse design studies that instantly find solutions have emerged \citep{liu2018nanophotonic, yilmaz2018airfoil, malkiel2018plasmonic}. 
The key characteristic of deep inverse design is that it focuses on generating optimal design candidates based on data, rather than searching for them from scratch.
This allows for real-time design optimization by implementing deep learning models on both simulation and optimization, given that the models are well-trained. 

While deep inverse design has demonstrated its potential for innovative design optimization through the studies in some specific areas \citep{liu2018nanophotonic, an2021multifunctionalgan}, there is still a lack of study that analyzes its general performance and feasibility.
Thus, this study aims to evaluate the performance of deep inverse design by comparing  with traditional design optimization.

This paper is structured as follows. 
In Section~\ref{sec:related-work}, the studies on traditional design optimization and deep inverse design are introduced. 
Section~\ref{sec:methods} proposes the methods for their comparison. 
Section~\ref{sec:perf-comparison} presents the resulting performance analysis and comparison.
The research findings are summarized and few directions for the field of deep inverse design are suggested in Section~\ref{sec:conclusion}.

The contributions of this study are as follow. 
First, we quantitatively evaluate the deep inverse design methods and compare them with the traditional design optimization methods. 
By performing an in-depth analysis on their performance, we anticipate to provide important indicators to establish effective and efficient design optimization processes for the field. 
Second, by analyzing the optimization performance with different problem definitions and data compositions, we suggest directions to enhance the usefulness of deep inverse design.

\section{Related Work}\label{sec:related-work}

\begin{figure*}[th]
    \centering
    \includegraphics[width=\textwidth]{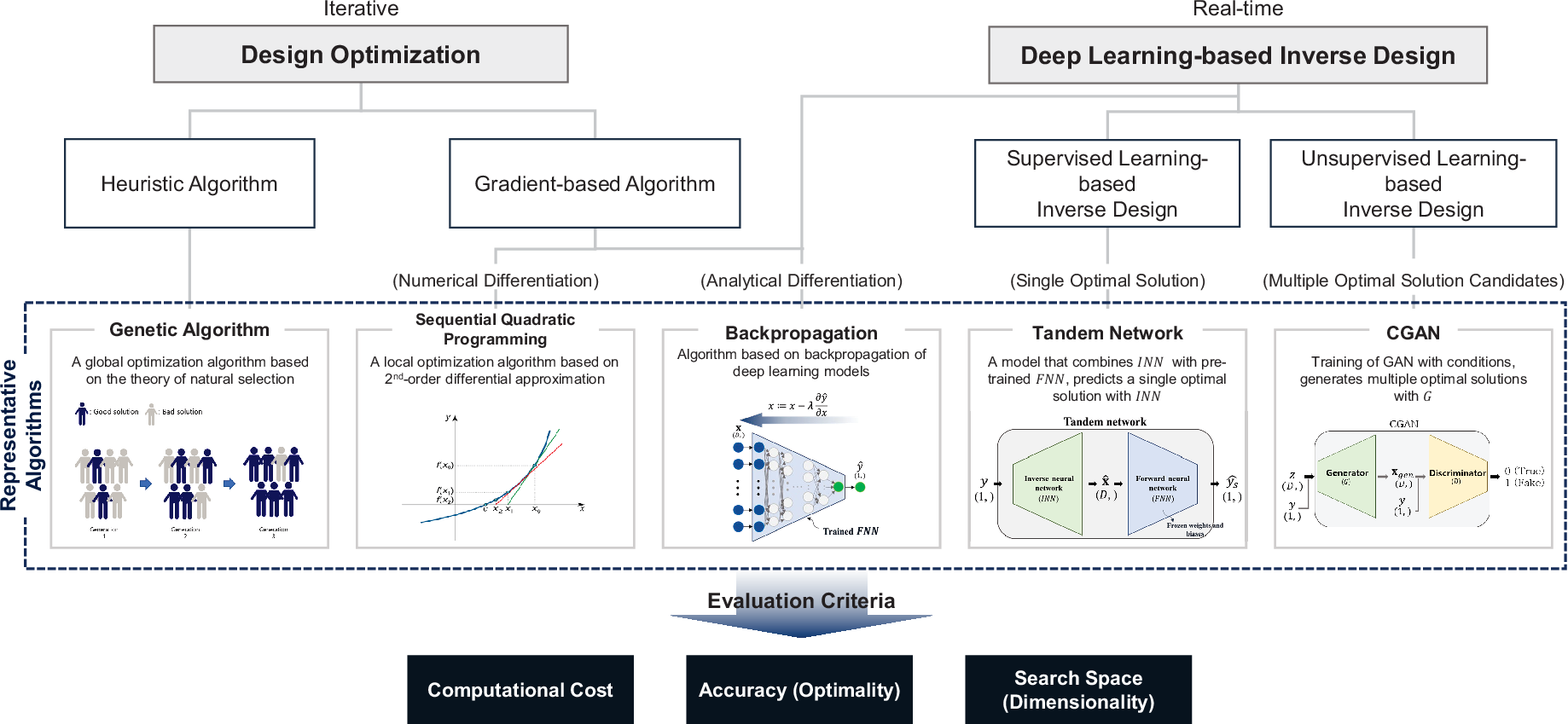}
    \caption{Classification of design optimization approaches and algorithms}
    \label{fig:classification-traditional-DO}
\end{figure*}

Many traditional design optimization and deep inverse design approaches can be classified as shown in Fig.~\ref{fig:classification-traditional-DO}. 
Traditional design optimization, which searches for optimal solutions iteratively, can be divided into heuristic and gradient-based algorithms. 
Deep inverse design, which predicts solutions in real-time, is divided into supervised and unsupervised learning based on how their neural networks are trained with data. 
Section~\ref{ssec:related-work-traditional-DO} introduces the representative models and applications of traditional design optimization, while \ref{ssec:related-work-supervised-deep-learning-inverse-design} and \ref{ssec:related-work-unsupervised-deep-learning-inverse-design} introduce the supervised and unsupervised learning-based inverse design methods, respectively.

\subsection{Traditional Design Optimization}\label{ssec:related-work-traditional-DO}

Traditional design optimization approaches can be divided into heuristic algorithms that generate and use random variables, and gradient-based algorithms that determine the search direction through differentiations. 
Genetic Algorithm (GA) \citep{whitley1994genetic} is one of the most representative models of heuristic approaches. 
It performs global optimization to generate a solution that satisfies an objective function by applying techniques such as cross-over or mutation based on the principle of natural selection. 
GA has been widely used in various practical fields such as transportation, architecture, and medicine to build optimal systems \citep{rahmani2011modeling, sharif2000multireservoir, nematzadeh2018medical}.
Sequential Quadratic Programming (SQP) \citep{gill2005snopt} and Backpropagation (BP) \citep{chen2020generative} are the representative models of gradient-based algorithms.
SQP performs local optimization by numerically approximating the second-order derivative of the objective function and constraints at each iteration.
On the contrary, BP performs local optimization by analytically computing the first-order derivative through the utilization of a surrogate neural network model.

These gradient-based algorithms have demonstrated their good optimization performance through many case studies in engineering design, especially in nanophotonics \citep{fesanghary2008hybridizing, peurifoy2017nanophotonic}.
Many studies have been conducted to systematically establish an understanding of optimization performance to use design optimization algorithms effectively \citep{aktemur2017comparison}. 
Further, a study combining global and local optimization algorithms had been proposed to pursue higher accuracy \citep{moussouni2007motor}.

\subsection{Supervised Deep Inverse Design}\label{ssec:related-work-supervised-deep-learning-inverse-design}

Supervised learning is a machine learning technique that learns the mapping relationship between input ($\mathbf{x}$) and output data ($y$). 
It is generally used to predict $y$ from $\mathbf{x}$, and depending on the type of output data, it is further classified into classification models and regression models.

A supervised learning model $f \colon \mathbf{x} \rightarrow y$ is trained to approximate a function, which outputs the predicted engineering performance ($\hat{y}$) for given design variables ($\mathbf{x}$) using a Forward Neural Network (FNN). 
This architecture has shown high prediction accuracies through many studies. 
In order to predict $\mathbf{x}$ that satisfy $y_t$, an Inverse Neural Network (INN) is used, which is a prediction model $f \colon y \rightarrow \mathbf{x}$.
However, this typically shows significantly lower prediction performance, as there may be multiple feasible design candidates ($\hat{\mathbf{x}})$ that correspond to a single $y_t$, failing to satisfy the one-to-one correspondence required for a traditional neural network.

To address this, \cite{liu2018nanophotonic} proposed a Tandem Network (TN) shown in Fig.~\ref{fig:tandem-network} that predicts an optimal design candidate satisfying the target spectral response. 
The proposed model combines a pre-trained FNN ($\mathbf{\hat{x}} \rightarrow \hat{y}$) and the INN ($y_t \rightarrow \mathbf{\hat{x}}$) that is trained thereafter.
By freezing the weights and biases of the pre-trained FNN that typically shows high prediction performance, TN is trained to solely optimize the parameters of the INN. 
As a result, $\mathbf{\hat{x}}$ of INN are passed to the frozen FNN which outputs $\hat{y}$, and INN is trained to minimize the loss between $y_t$ and $\hat{y}$.
This enables the INN to generate an optimal $\mathbf{\hat{x}}$ for $y_t$. 

TN has been mainly used in nanophotonics to design nano-scale optical structures. 
Optimizing fine structures typically involves solving highly complex non-convex problems.
Although global optimization methods have been studied for nanophotonics, there have been limitations in handling excessive computational costs \citep{molesky2018nanophotonics}. 
To address this, the studies applying deep inverse design have emerged, predicting design parameters such as the thickness of nanostructures that can lead to target wavelengths or spectral responses \citep{peurifoy2018nanophotonic, malkiel2021nanostructures}.
In addition, the studies have been proposed to design the shape of disc brakes that satisfy target APT and drag performance \citep{kim2022brakes}, and the shape of an airfoil that satisfies the target lift \citep{sekar2019airfoil}.

\subsection{Unsupervised Deep Inverse Design}\label{ssec:related-work-unsupervised-deep-learning-inverse-design}

Unsupervised learning is a machine learning technique that learns patterns and characteristics of unlabeled data. 
It is commonly used for tasks such as clustering or dimensionality reduction, and also for generative modeling.
One of the most representative unsupervised deep inverse design models is Generative Adversarial Networks (GAN), proposed by \cite{goodfellow2020generative}, which is a model that generates new data by learning the characteristics and distribution of the data.

GAN is trained  sequentially by optimizing the parameters of the generator network ($G$) and the discriminator network ($D$).
GAN has evolved in various fields such as cancer prognosis prediction and speech synthesis \citep{zhu2020cancer, kumar2019melgan}.
Many GAN variants have been developed for different purposes, such as Conditional GAN (CGAN) \citep{mirza2014conditional}.

\begin{equation}\label{eq:CGAN}
\begin{aligned}
    \min_{G} \max_{D} V(D,G) \\&= \mathbb{E}_{x \sim p_{data}(\mathbf{x})} \left[ \log D(\mathbf{x} | y) \right] \\&+ \mathbb{E}_{z \sim p_z (z)} \left[ \log (1 - D(G(z|y))) \right]
\end{aligned}        
\end{equation}

CGAN additionally incorporates conditional information by giving them as the inputs to $G$. 
It is trained using the modified objective function defined in Eq.~\eqref{eq:CGAN} with condition terms ($y$) added, which is commonly used to generate data that belong to specific classes or to satisfy a target performance.

The aim of unsupervised inverse design is to incorporate the conditions in order to learn the relationship between design variables and engineering performance.
The model is capable of generating many unique design candidates based on the conditional input of target performance. 
As a result, CGAN has been actively used for inverse design purposes, such as designing airfoil shapes \citep{yilmaz2020airfoilcgan} and meta-materials \citep{an2021multifunctionalgan}. 
To improve the accuracy of the model, the studies proposing modifications on CGAN that incorporate different ways of training have been conducted extensively \citep{dai2022cgan, mall2020cyclical, dong2020graphenecgan}.

\section{Methods}\label{sec:methods}

\begin{figure*}[tb]
    \centering
    \includegraphics[width=\textwidth]{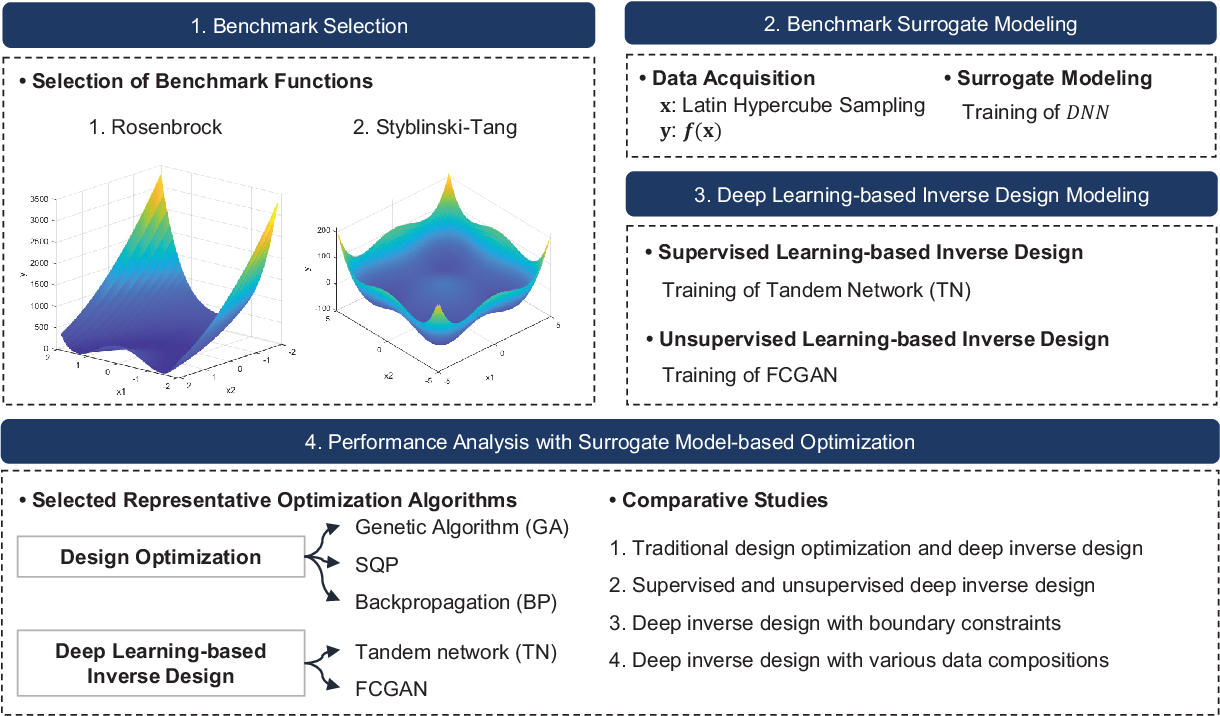}
    \caption{Research framework}
    \label{fig:research-framework}
\end{figure*}

In order to analyze deep inverse design and compare with traditional design optimization, this study proposes a framework shown in Fig.~\ref{fig:research-framework} that builds deep inverse design models and analyzes their performance for surrogate model-based optimization problems based on the benchmark functions.
A brief explanation of each stage is provided below.

\begin{enumerate}
    \item Benchmark selection: 
    the two multi-dimensional non-convex functions, Rosenbrock and Styblinski-Tang functions, are selected for the benchmark functions to evaluate the traditional optimization and deep inverse design models.
    \item Benchmark surrogate modeling: 
    the surrogate models for each benchmark were constructed. 
    The models are based on a neural network trained with the data obtained from the benchmark functions. 
    \item Deep learning-based inverse design modeling: 
    TN and FCGAN were selected and trained for the representative supervised and unsupervised deep inverse design models, respectively.
    \item Performance analysis with surrogate model-based optimization: 
    first, numerous surrogate model-based optimization problems are defined to evaluate the optimization performance of the selected models.
    To analyze their performance from multiple perspectives, four comparative studies are designed.
\end{enumerate}

\subsection{Benchmark Selection}\label{ssec:benchmark-selection}
The aim of this study is not to build accurate surrogate models. 
Instead, we aim to compare the performance of the traditional design optimization and deep inverse design models when performing surrogate model-based optimization. 
Therefore, the models are optimized on the same surrogate models, and their performance is evaluated.
Furthermore, since the performance of the models depends on the search space, it is necessary to design various environments for a comparison.

To evaluate the performance of a new optimization method, it is necessary to compare and analyze its performance against existing optimization methods with benchmark problems with various characteristics and dimensions. 
\cite{li2013benchmark} and \cite{molga2005testfunction} proposed few mathematical evaluation functions that can be used to evaluate the optimization performance. 
These proposed functions range from simple convex shapes (e.g., Sphere function) to complex shapes with multiple valleys (e.g., Rastrigin function).
These shapes are typically represented with different characteristics such as dimensionality, modality, and separability. 
Below is the explanation of each key characteristic \citep{jamil2013literature}:
\begin{enumerate}
    \item Dimensionality: 
    the number of design variables. 
    As the number of design variables increases, the search space grows exponentially, making optimization more challenging.
    \item Modality: 
    the number of peaks in the shape of a function. 
    The more peaks there are, the higher the likelihood that optimization algorithms will converge to local optima.
    \item Separability (decomposability): 
    the independence of design variables from each other. 
    Generally, optimization becomes more difficult when the design variables have high dependencies. 
\end{enumerate}

The mathematical functions that follow the diverse characteristics of real-world engineering design problems must be selected as the benchmarks for this study. 
To ensure the difficulty of finding global optimal solutions and to incorporate different levels of complexity, we considered two non-convex functions with multiple modes. 
Since the size of search space depends on the domain, we also considered functions that allow dimensionality expansion to compare performance from low to high dimensions. 
As a result, the selected benchmark functions are Rosenbrock function and Styblinski-Tang function. 
Table~\ref{tab:benchmark-functions} shows their characteristics.
The detailed explanations of the two benchmark functions are as follow.

\begin{table}[tb]
\caption{The characteristics of the selected benchmark functions}\label{tab:benchmark-functions}%
\begin{tabular}{@{}ccc@{}}
\toprule
function & Rosenbrock & Styblinski-Tang \\
\midrule
modality & single & multi \\
separability & no & yes \\
convexity & \multicolumn{2}{c}{non-convex} \\
dimensionality & \multicolumn{2}{c}{multi} \\
\botrule
\end{tabular}
\end{table}

\paragraph{Rosenbrock Function}

\begin{figure}[htb]
    \centering
    \includegraphics[width=0.4\textwidth]{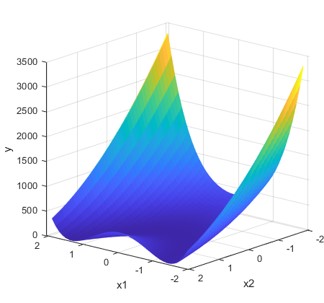}
    \caption{Rosenbrock function in 2D}
    \label{fig:rosenbrock-shape}
\end{figure}

\begin{equation}\label{eq:rosenbrock-function}
    f_{\mathrm{Rosenbrock}}(x) = \sum^{D-1}_{i=1} \left[ 100(x_{i+1} - x_{i}^2)^2 + (1 - x_i)^2 \right]
\end{equation}

Rosenbrock function defined in Eq.~\eqref{eq:rosenbrock-function} is a non-convex function for basic optimization benchmark purposes, also known as the second function of De Jong \citep{rosenbrock1960automatic}. 
It has the shape of a narrow, elongated parabolic valley as shown in Fig.~\ref{fig:rosenbrock-shape}.
The design space is usually limited to $x \in [-2.0,...,2.0]^D$, where $D$ refers to the dimensionality of the data.

\paragraph{Styblinski-Tang Function}

\begin{figure}[htb]
    \centering
    \includegraphics[width=0.4\textwidth]{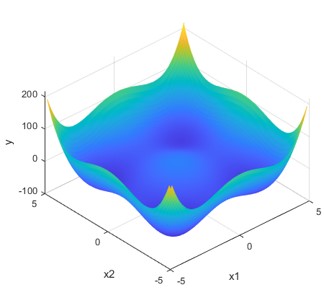}
    \caption{Styblinski-Tang Function in 2D}
    \label{fig:styblinski-tang-shape}
\end{figure}

\begin{equation}\label{eq:styblinski-tang-function}
    f_{\mathrm{Styblinski-Tang}}(x) = \frac{1}{2} \sum^{D}_{i=1} (x_i^4 - 16x_i^2 + 5x_i)
\end{equation}

Styblinski-Tang function defined in Eq.~\eqref{eq:styblinski-tang-function} is the evaluation function proposed by \cite{szu1986tang} and modified by \cite{styblinski1990nonconvex}.
It has a non-convex shape with four local minima, as shown in Fig.~\ref{fig:styblinski-tang-shape}.
The design space is often limited to $x \in [-5.0,...,5.0]^D$. 

Both of the selected functions are multi-dimensional. 
However, while Rosenbrock function is uni-modal and non-separable, Styblinski-Tang function is multi-modal and separable.
The dimensionality defines the search space of the design variables, which affects the optimization performance. 
Thus, the optimization performance of the models will be analyzed with the benchmark functions from low to high dimensions. 
Six different cases of dimensionality were defined for each benchmark: 2, 10, 30, 50, 70, and 100 dimensions.

\subsection{Benchmark Surrogate Modeling}\label{ssec:benchmark-surrogate-modeling}
Surrogate model-based optimization replaces time-consuming simulations in the design optimization process with a predictive surrogate model.
Thus, it is often preferred in the field due to its time efficiency. 
Below explains the deep learning-based surrogate modeling process for the benchmark functions.

\subsubsection{Data Acquisition and Preprocessing}
The data from the two selected benchmark functions are needed to train the deep learning-based surrogate models. 
First, Latin Hypercube Sampling (LHS) was used to sample the training data.
LHS is a statistical method that allows to sample over the search space of a multi-dimensional distribution \citep{iman1980lhs}. 
For each function and dimension, we sample 10,000 data points ($\mathbf{x}$) that satisfy the boundary constraints. 
and $y$ that correspond to $\mathbf{x}$ were collected using the ground truth benchmark functions. 
The resulting $\mathbf{x} \colon y$ pair data can be summarized as shown in Table~\ref{tab:benchmark-data-info-2}.

\begin{table}[tb]
\caption{The characteristics of the data from the benchmark functions used to train the surrogate models}\label{tab:benchmark-data-info-2}%
\begin{tabular}{@{}ccc@{}}
\toprule
function & Rosenbrock & Styblinski-Tang \\
\midrule
boundary & $[-2.0, ... ,2.0]^D$ & $[-5.0,...,5.0]^D$ \\
dimension & \multicolumn{2}{c}{2, 10, 30, 50, 70, 100}  \\
data size & \multicolumn{2}{c}{10,000} \\
\botrule
\end{tabular}
\end{table}

The preprocessing of the generated data is as follows. 
First, the data normalization was performed. 
$\mathbf{x}$ was scaled to values between 0 and 1 using min-max scaling, and $y$ was normalized to have a distribution with mean and variance of 0 and 1, respectively. 
Next, a data split was performed on 10,000 data points, building a train dataset (9,000) and a validation dataset (1,000) with a 9:1 ratio for model training and hyperparameter optimization, respectively.


\subsubsection{Deep Learning-based Surrogate Modeling}\label{ssec:deep-surrogate-modeling}

\begin{figure}[tb]
    \centering
    \includegraphics[width=0.4\textwidth]{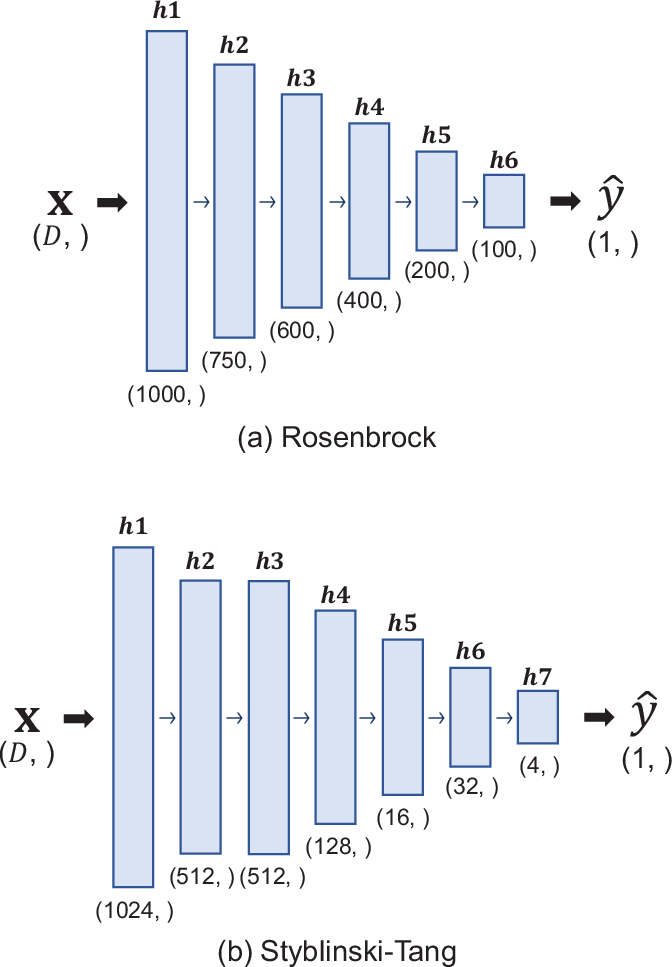}
    \caption{The architectures of the deep surrogate models for the benchmark functions}
    \label{fig:surrogate-architecture}
\end{figure}

The surrogate models for each benchmark function were built by training a deep neural network with the generated dataset.
The parameters of each model were optimized by minimizing the Mean Squared Error (MSE)  shown in Eq.~\eqref{eq:MSE}.

\begin{equation}\label{eq:MSE}
    \mathrm{MSE} = \frac{1}{n} \sum^{n}_{i=1} (y_i - \hat{y}_i)^2
\end{equation}
$y_i$ and $\hat{y}_i$ are the ground truth and the estimation by a model, respectively.
The network architecture and hyperparameter optimization were performed to improve the prediction performance. 
Fig.~\ref{fig:surrogate-architecture} shows the final network architecture.

To prevent overfitting the surrogate models, early stopping was applied and the dropout layers were added to the model architecture. 
Table~\ref{tab:surrogate-model-hyperparameters} shows the final hyperparameters used to train the models.
The surrogate models were constructed for six different dimensions per benchmark, resulting in a total of twelve models.

The designed surrogate models showed good performance while training.
The authors refer to Appendix~\ref{app:surrogate-models} for the training results of the surrogate models.

\begin{table}[tb]
\caption{The hyperparameters of the surrogate models}\label{tab:surrogate-model-hyperparameters}%
\begin{tabular}{@{}cc@{}}
\toprule
hyperparameters & values\\
\midrule
learning rate & 0.0001 \\
epoch & 1,000 \\
batch size & 128 \\
optimizer & ADAM \\
activation function & ELU or ReLU \\
\botrule
\end{tabular}
\end{table}

\subsubsection{Surrogate Model-based Optimization}\label{ssec:methods-surrogate-evaluation}

The evaluation criteria were chosen for the performance of the traditional design optimization and deep inverse design models with surrogate model-based optimization. 
An optimization problem is typically defined using design variables, objective functions and constraints. This section demonstrates these key elements.

\paragraph{Design variables}
\begin{equation}
    \mathbf{x} \subset \mathbb{R}^D    
\end{equation}

The search space of design variables is defined as a set of real numbers in the dimension ($D$).

\paragraph{Objective function}
\begin{equation}\label{eq:objective-function}
    \min_\mathbf{x} \mathrm{RMSE} (y_t, \hat{y}(\mathbf{x}))
\end{equation}

The objective function is set to minimize the squared root of MSE between $y_t$ and $\hat{y}$ as shown in Eq.~\eqref{eq:objective-function}, where $y_t$ represents the target value, and $\hat{y}$ represents the prediction for $\mathbf{x}$ through the surrogate models. 
Thus, the models must aim to find $\mathbf{x}$ that makes $\hat{y}$ closer to $y_t$.

\paragraph{Constraints}
\begin{equation}
    \mathbf{x} \in \left[ lb, ub \right]^D
\end{equation}

The constraints are the boundary limits of the design variables, where $lb$ and $ub$ are the lower and upper limits. 
The constraints are defined to enable the exploration of optimal solutions within these ranges.
Table~\ref{tab:benchmark-data-info-2} shows the boundary limits for the benchmark functions.
These limits will also be used to see whether the models are able to provide $\hat{\mathbf{x}}$ that satisfy these constraints.

\paragraph{Evaluation Criteria}

\begin{equation}\label{eq:acc}
    \mathcal{L} = \frac{1}{n} \sum^{n}_{i=1} \mathrm{RMSE}(y_{i,t}, \hat{y}(\mathbf{x}_i^*))
\end{equation}

Eq.~\eqref{eq:acc} shows the loss function ($\mathcal{L}$) to assess the model performance, where $y_{i,t}$ is the target performance and $\mathbf{x}_i^*$ is the corresponding design predictions.

The optimized models will be evaluated on the test dataset consisting of 100 data points ($n=100$) directly sampled from the benchmark functions using LHS for each dimension.

The performance of the traditional design optimization methods typically depends on the hyperparameters for training, such as initial solutions or the maximum number of iterations. 
To ensure a fair comparison between models, the same optimization parameters were used for GA, SQP and BP to minimize their impacts as shown in Table~\ref{tab:optimization-params}.
The initial points were generated using LHS to ensure that the optimal solution could be explored from similar conditions.

\begin{table}[tbh]
\caption{The parameters for optimization}\label{tab:optimization-params}%
\begin{tabular}{@{}cccc@{}}
\toprule
optimization parameters & GA & SQP & BP\\
\midrule
Num. of initial solution & \multicolumn{3}{c}{10}\\
Max. iterations & \multicolumn{3}{c}{100}\\
\botrule
\end{tabular}
\end{table}

\subsection{Deep Inverse Design Modeling}\label{ssec:deep-inverse-design-modeling}
\subsubsection{Modeling of Deep Inverse Design}\label{sssec:modeling-of-deep-learning-based-inverse-design}

TN and FCGAN were selected to represent the supervised and unsupervised inverse design models, respectively.
Below explains their structures and training process.

\paragraph{Tandem Network (supervised)}

\begin{figure}[tb]
    \centering
    \includegraphics[width=0.5\textwidth]{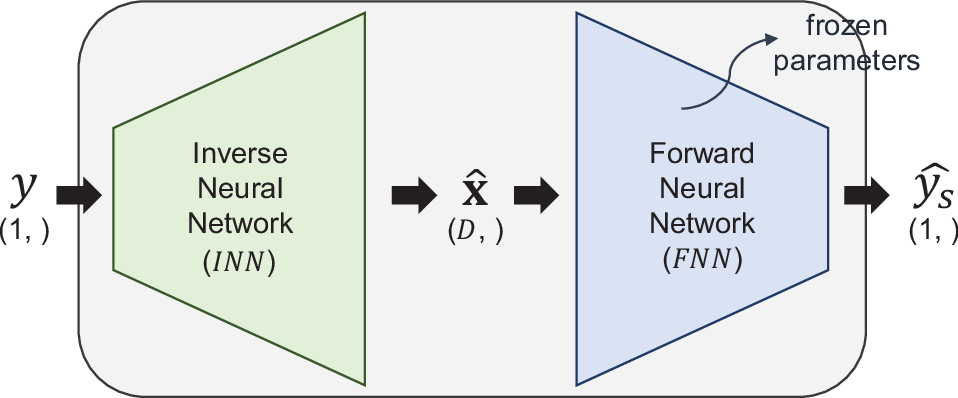}
    \caption{The structure of Tandem Network}
    \label{fig:tandem-network}
\end{figure}

TN has a structure that combines a pre-trained Forward Neural Network (FNN) and an Inverse Neural Network (INN) as shown in Fig.~\ref{fig:tandem-network}.
In this study, the FNN of TN is the trained surrogate model.
The network structures of INN is a mirrored version of FNN for each benchmark function shown in Fig.~\ref{fig:surrogate-architecture}.
The parameters of the FNN from the pre-trained surrogate model are frozen, so only the parameters of the INN are optimized during training.

\paragraph{FCGAN (unsupervised)}

\begin{figure}[tb]
    \centering
    \includegraphics[width=0.5\textwidth]{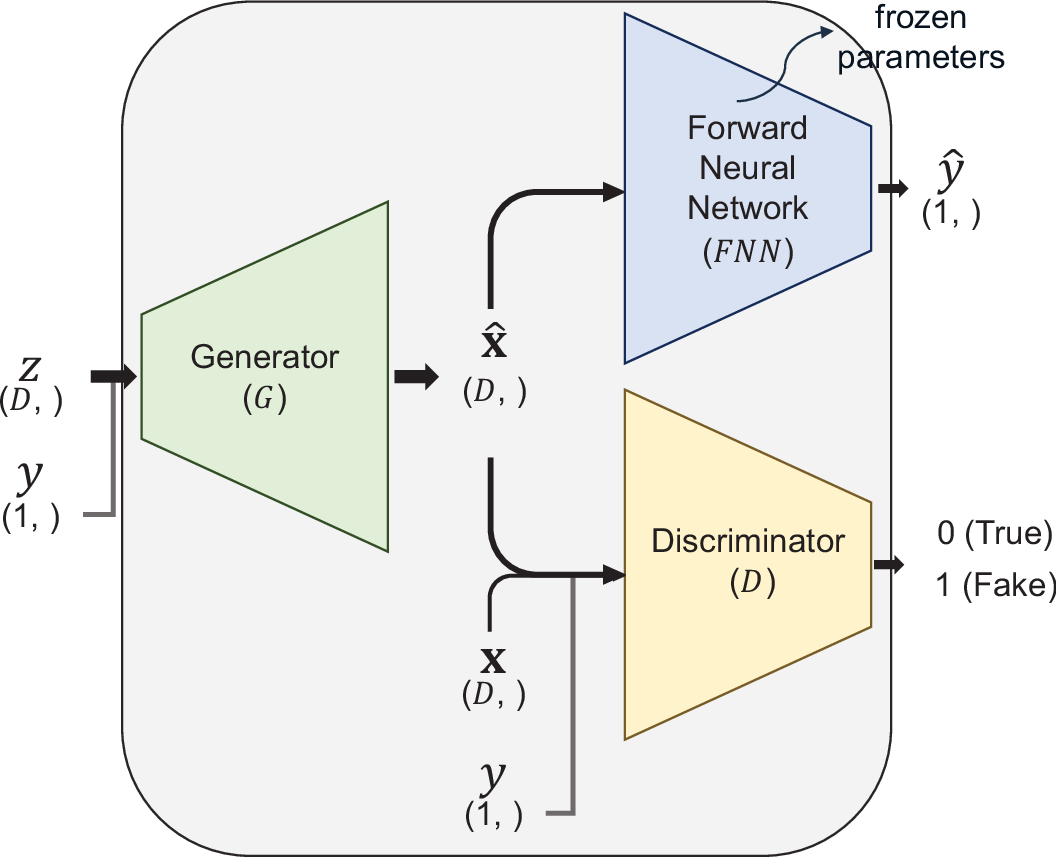}
    \caption{The structure of FCGAN}
    \label{fig:FCGAN}
\end{figure}

Forward and Conditional Generative Adversarial Network (FCGAN) combines FNN ($F$) with the baseline CGAN as shown in Fig.~\ref{fig:FCGAN}.
Similar to TN, $F$ in the FCGAN is also a pre-trained surrogate model with its parameters frozen during the training.
CGAN simultaneously trains the generator ($G$), which aims to optimize Eq.~\eqref{eq:CGAN} by fooling the discriminator ($D$).
In FCGAN, $G$ is trained to minimize the loss between $y_t$ and the $\hat{y}$ estimated by $F$ from $\hat{\mathrm{x}}$ that $G$ generates.
Thus, $G$ is trained to generate $\hat{\mathrm{x}}$ that, when passed through FNN, produces $\hat{y}$ similar to the conditional target input ($y_t$). 
Eq.~\eqref{eq:FCGAN} is the resulting multi-objective loss function of FCGAN.
The structures of $G$ and $D$ are identical to the structures of INN and FNN of TN.

\begin{equation}
    \label{eq:FCGAN}
    \begin{aligned}
        \min_{G} \max_{D} V(D,G) \\&= \mathbb{E}_{\mathrm{x} \sim p_{data}(\mathrm{x})} \left[ \log D(\mathrm{x} | y) \right] \\&+ \mathbb{E}_{z \sim p_z (z)} \left[ \log (1 - D(G(z|y))) \right] \\&+ \mathrm{MSE}(y, F(G(z|y)))
    \end{aligned}
\end{equation}

\subsubsection{Modeling with Boundary Constraints}\label{sssec:methods-bc}
Most engineering design problems involve boundary constraints. 
Traditional design optimization models are capable of considering these boundary constraints when looking for optimal solutions. 
However, conventional deep learning-based models rely solely on data when training and the search space is not explicitly specified, which may lead to providing solutions outside boundaries. 
To address this issue, a deep learning architecture is proposed that considers boundary constraints by applying an activation function on its output, allowing for the prediction of solutions within a specified range. 

\begin{figure}[tb]
    \centering
    \includegraphics[width=0.4\textwidth]{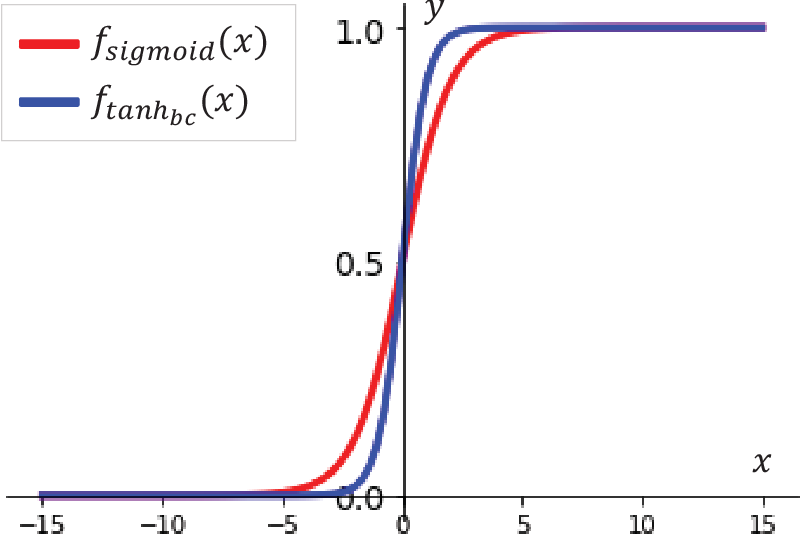}
    \caption{sigmoid and $\tanh_{bc}$ activation functions to limit the prediction range of the deep inverse design models}
    \label{fig:tanh_sigmoid_vis}
\end{figure}

During the training of the surrogate models, $\mathbf{x}$ was scaled to a range between 0 and 1 using min-max scaling. 
By using an activation function that outputs values between 0 and 1, the range of predictions made by a model can be restricted. 
The activation functions used are the sigmoid function in Eq.~\eqref{eq:sigmoid} and a modified version of the tanh function, 
referred to as $\tanh_{bc}$ function in Eq.~\eqref{eq:tanh_bc}. 
These activation functions enable the prediction of design variables within the range [0, 1] as shown in Fig.~\ref{fig:tanh_sigmoid_vis}.
They were applied to the output layers of TN and FCGAN to create deep inverse design models with boundary constraints, referred as TN-BC and FCGAN-BC in this paper.

\begin{equation}\label{eq:sigmoid}
    f_{\mathrm{sigmoid}}(x) = \frac{1}{1 + e^{-x}}
\end{equation}

\begin{equation}\label{eq:tanh_bc}
    f_{{\tanh}_{bc}}(x) = \frac{1}{2} \left( \frac{e^{x} - e^{-x}}{e^{x} + e^{-x} + 1} \right)
\end{equation}

The designed deep inverse design models showed good performance while training.
The authors refer to Appendix~\ref{app:deep-models} for the training results of the models.

\subsection{Case Studies for Comparison}\label{ssec:methods-comparative-studies}
This study aims to analyze the optimization performance of design optimization and deep inverse design from many perspectives. 
To achieve this aim, four case studies were designed as below.

\paragraph{Case 1: design optimization and deep inverse design}
The deep inverse design models are compared to the traditional design optimization models which have been widely studied and used in the field. 
The most significant difference between them is how optimal design solutions are obtained. 
The robust performance of traditional methods has been validated through numerous studies. 
However, they typically require high computational costs as they iteratively search for solutions.
In contrast, deep inverse design models are capable of instantly predicting optimal solutions after the training, but there have been lacking studies on validating their performance.
Therefore, this comparative study must be conducted to evaluate their optimization performance in terms of accuracy and computational costs.

\paragraph{Case 2: supervised and unsupervised deep inverse design}
Most deep inverse design models can be classified as supervised or unsupervised learning-based. 
The main difference between these two is the inclusion of $y$ for their training, although unsupervised models may also consider $y$ as conditional inputs ($y_t$) when predicting $\hat{\mathbf{x}}$.
While supervised inverse design models focus only on the relationship between $\mathbf{x}$ and $y$, unsupervised inverse design models learn both the relationship between $\mathbf{x}$ and $y$ and the distribution of $\mathbf{x}$ simultaneously. 
Thus, this case study aims to observe how the accuracy and diversity of an unsupervised inverse design model differ from a supervised model. 
Furthermore, various weight combinations were applied to FCGAN that uses multi-objective functions and the diversity and accuracy of $\hat{\mathbf{x}}$ were measured for each combination.

\paragraph{Case 3: deep inverse design with boundary constraints}
TN-BC that considers boundary constraints using an activation function was proposed in Section~\ref{sssec:methods-bc}. 
The performance of TN-BC was analyzed and compared to that of the baseline TN on whether their outputs satisfy the boundary constraints. 
Through this case study, the applicability of deep inverse design can be inferred to real-world design problems where boundary constraints are often essential.

\paragraph{Case 4: deep inverse design with various data compositions}
In the real world, it is often challenging to collect a large amount of data for surrogate modeling. 
Thus, it is common to perform over-sampling from a small number of real observation data samples to augment the dataset. 
Since deep inverse design models rely on data for training, their performance can be hugely affected by the size and ratio of the real and augmented data. 
Therefore, this case study aims to compare of the optimization performance of the models with various data compositions.

\section{Performance Comparison}\label{sec:perf-comparison}
In this section, we present the results of the four case studies introduced in Section~\ref{ssec:methods-comparative-studies}.
The authors refer to Appendix~\ref{app:model-training-results} for the training results of the deep surrogate models and the deep inverse design models.

\subsection{Traditional Design Optimization and Deep Inverse Design}

This section compares the optimization performance of the design optimization and deep inverse design methods in various dimensions.
The selected criteria for comparing their performance are the accuracy and computational cost.
The accuracy of the models is evaluated using average RMSE values calculated using Eq.~\eqref{eq:acc} for 100 test data points ($n=100$). 
The computational cost in this study is defined as the number of model algorithm iterations. 
For GA, it is the number of generations until the optimal solution is obtained.
For SQP and BP, it is the number of iterations until the optimal solution is achieved. 
For deep inverse design models, it is the number of predictions.

\begin{figure}[tb]
    \centering
    \includegraphics[width=0.5\textwidth]{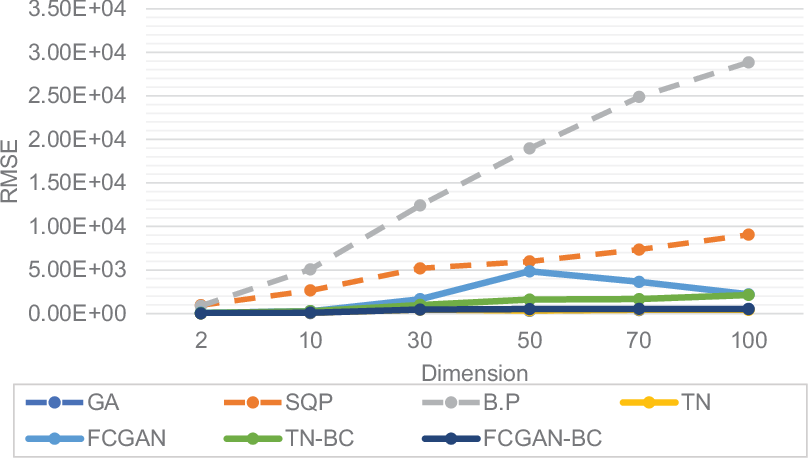}
    \caption{The RMSEs of the models for Rosenbrock function as the dimension increases}
    \label{fig:RMSEs-Rosenbrock-dim-29}
\end{figure}

\begin{figure}[tb]
    \centering
    \includegraphics[width=0.5\textwidth]{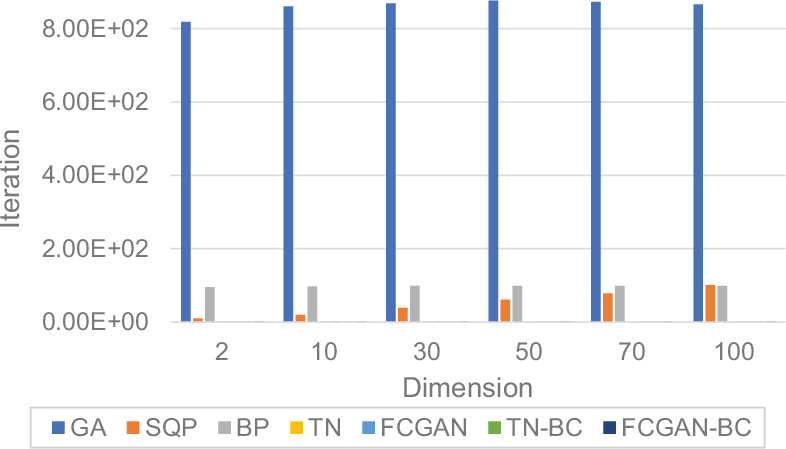}
    \caption{The computational costs of the models for Rosenbrock function as the dimension increases}
    \label{fig:comp-costs-Rosenbrock-dim-30}
\end{figure}

\begin{figure}[tb]
    \centering
    \includegraphics[width=0.5\textwidth]{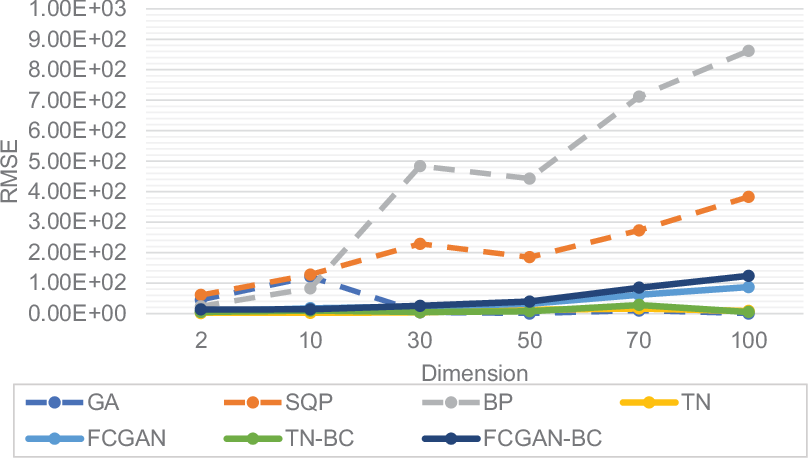}
    \caption{The RMSEs of the models for Styblinski-Tang function as the dimension increases}
    \label{fig:RMSEs-Styblinski-dim-31}
\end{figure}

\begin{figure}[tb]
    \centering
    \includegraphics[width=0.5\textwidth]{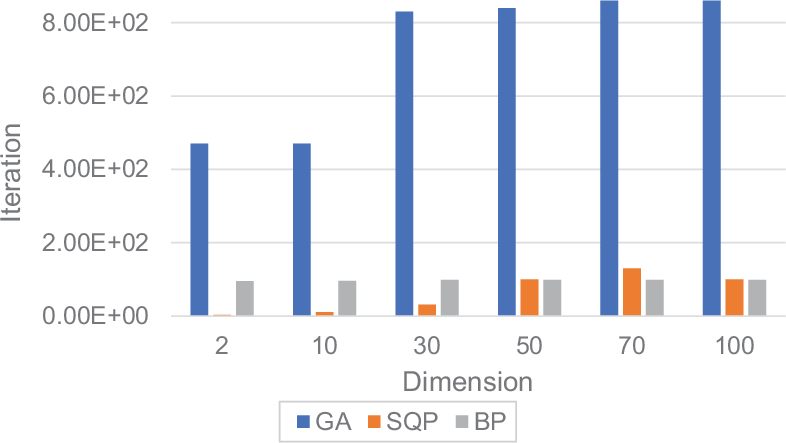}
    \caption{The computational costs of the models for Styblinski-Tang function as the dimension increases}
    \label{fig:comp-costs-Styblinski-dim-32}
\end{figure}

The accuracies and the computational costs of the models for Rosenbrock function are plotted in Fig.~\ref{fig:RMSEs-Rosenbrock-dim-29} and \ref{fig:comp-costs-Rosenbrock-dim-30}.
The accuracies and the computational costs of the models for Styblinski-Tang function are plotted in Fig.~\ref{fig:RMSEs-Styblinski-dim-31} and \ref{fig:comp-costs-Styblinski-dim-32}.
Note that the computational costs for the deep learning models are too low to be visible in Fig.~\ref{fig:comp-costs-Rosenbrock-dim-30} and \ref{fig:comp-costs-Styblinski-dim-32} as they only take a single iteration for inference.

As the dimensionality increases, the size of the search space also increases, typically leading to performance degradation in terms of both accuracy and computational cost. 
The performance and characteristics of the models across multiple dimensions are analyzed below.

First, we analyze the accuracies of the models. 
GA showed high accuracy in all dimensions. 
In particular, for the complex 100D Styblinski-Tang function, its RMSE was near zero for the test dataset. 
In contrast, the local optimization algorithms, SQP and BP, showed the lowest accuracies.
They suffered a significant performance degradation as the search space increased. 
Furthermore, the deep unsupervised model, FCGAN, showed a gradual performance degradation with as the dimension increases, while the supervised TN showed very high accuracy similar to GA. 
The unsupervised model showed relatively weaker performance compared to the other models in large and complex search spaces.

Second, we analyze the computational costs of the models.
GA, which demonstrated high accuracy, consumed significantly high computational cost in all dimensions. 
SQP and BP consumed a moderate computational cost. 
SQP consumes more computational cost as the dimensionality increased to achieve convergence, while BP was not affected by the dimensionality as its optimization only involves gradient calculations for a fixed number of variables. 
The deep learning-based models consumed the least regardless of dimensionality or function since it can derive the optimal solution with a single iteration.

In general, the traditional optimization models are significantly affected by dimensionality, whereas the deep inverse design models demonstrated good accuracies and small computational costs regardless of dimensions. 
In particular, they showed good performance in predicting solutions in real-time for Styblinski-Tang function with complex search spaces.

\begin{figure}[tb]
    \centering
    \includegraphics[width=0.5\textwidth]{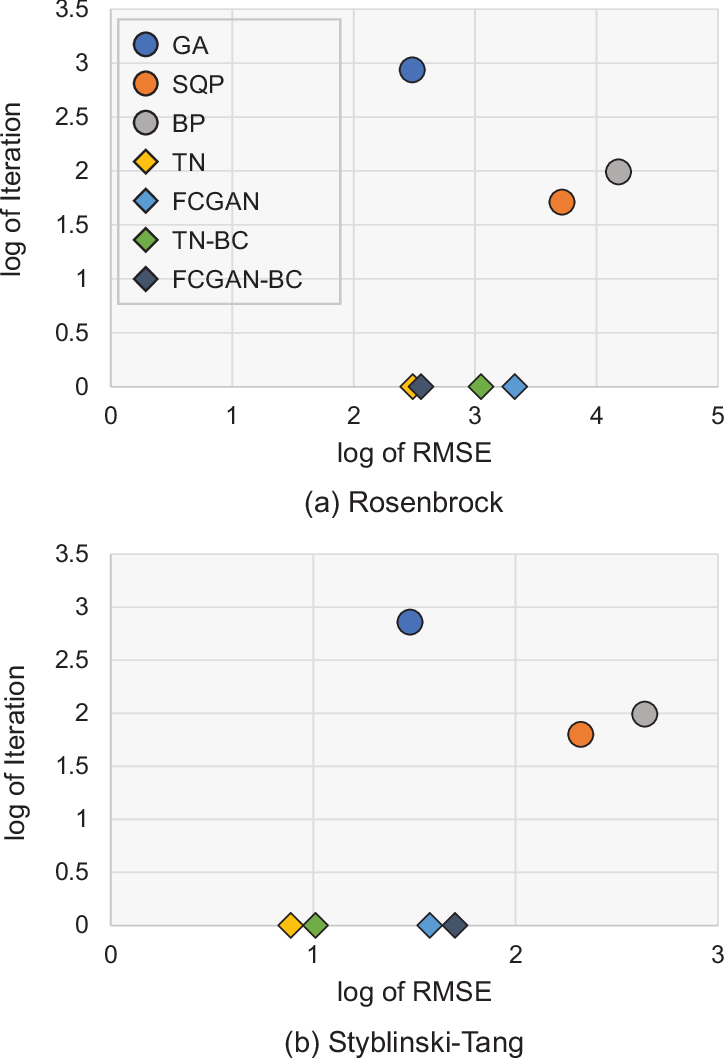}
    \caption{The logarithms of RMSEs and the computational costs of the models for the benchmark functions. The lower the values, the better the performance.}
    \label{fig:scatter-opt-perf-33-34}
\end{figure}

Fig.~\ref{fig:scatter-opt-perf-33-34} shows the scatter plots of the RMSE values against the computational costs averaged across the dimensions for each benchmark function.
GA may be considered the best model in terms of accuracy, but it required very high computational cost to achieve high accuracy. 
SQP and BP, the gradient-based algorithms, consumed moderate costs for both benchmarks but showed the worst levels of accuracy. 
The deep learning-based models demonstrated good performance in terms of accuracy and computational cost for both benchmark functions.

The authors refer to the tables in Appendix~\ref{app:perf-comparison} for the values used to plot the figures in this section.

\subsection{Supervised and Unsupervised Inverse Design}
\subsubsection{Performance Comparison}
Supervised learning learns from labeled data which enables the derivation of accurate solutions, whereas unsupervised learning is capable of deriving diverse solutions for a single target.
In this case study, the accuracy and diversity were chosen as the evaluation criteria to compare the performance of the models. 
Fig.~\ref{fig:scatter-opt-perf-33-34} demonstrates the superior accuracy of the deep inverse design models for the simple-shaped Rosenbrock function. 
However, Fig.~\ref{fig:RMSEs-Styblinski-dim-31} and \ref{fig:scatter-opt-perf-33-34} show that the performance of the unsupervised FCGAN for Styblinski-Tang function, which do not heavily rely on labeled data, is lower than that of the supervised TN. 
This emphasizes the importance of learning from labeled data in complex optimization problems.

\begin{figure}[tb]
    \centering
    \includegraphics[width=0.4\textwidth]{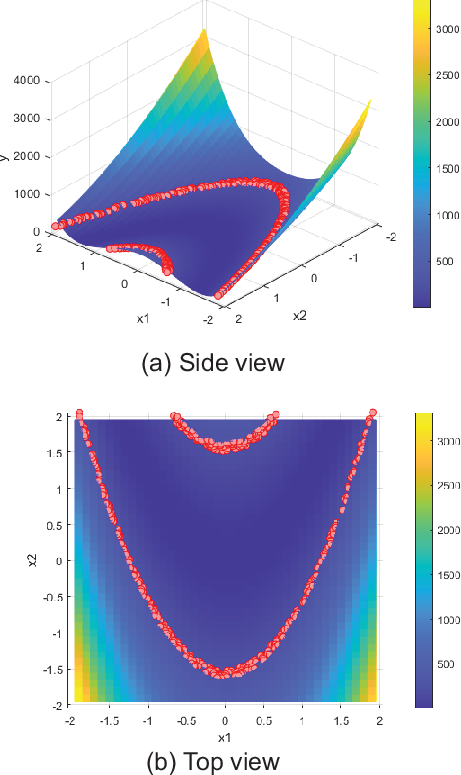}
    \caption{The distribution of sampled $\mathbf{x}$ for $y_t$ in Rosenbrock function}
    \label{fig:dist-x-35}
\end{figure}

\begin{figure}[tb]
    \centering
    \includegraphics[width=0.5\textwidth]{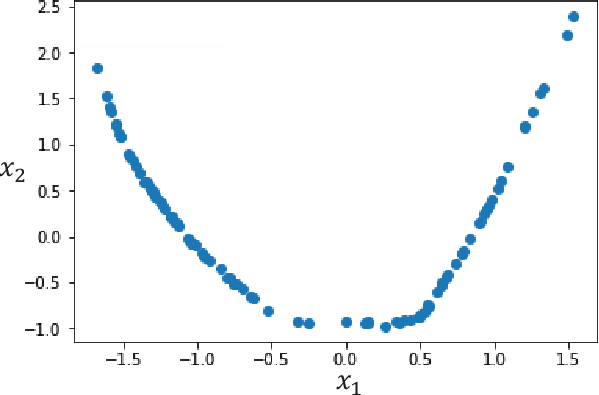}
    \caption{The distribution of $\hat{\mathbf{x}}$ generated by FCGAN}
    \label{fig:scatter-FCGAN-x-36}
\end{figure}

While unsupervised learning-based models may exhibit lower accuracy compared to supervised learning-based techniques, they have a significant advantage in generating multiple new predictions ($\hat{\mathbf{x}}$). 
To illustrate this, $\hat{\mathbf{x}}$ were sampled from the FCGAN model for Rosenbrock function in 2D, generated to predict the specific target ($y_t$). 
$y_t$ was set to 250, while the ground truth ($\mathbf{x}$) was set as the distribution that satisfy this target within an error range of $\pm 20$. 
Fig.~\ref{fig:dist-x-35} shows the distribution of $\mathbf{x}$ with the range marked as red. 
By training the FCGAN model with the randomly sampled noise vectors ($z$) and the conditional input ($y_t$), 100 $\hat{\mathbf{x}}$ values were generated. 
Fig.~\ref{fig:scatter-FCGAN-x-36} shows the distribution of $\hat{\mathbf{x}}$ as a scatter plot.
It can be seen that the distribution of $\hat{\mathbf{x}}$ is similar to that of $\mathbf{x}$ in Fig.~\ref{fig:dist-x-35}.

\subsubsection{Weights for Multi-objective Functions}
FCGAN is trained through the multi-objective function optimization of $D$ and FNN. 
As $D$ is trained, the learning capability of $G$, which generates $\hat{\mathbf{x}}$, also improves, allowing the model to better learn the distribution of $\mathbf{x}$. 

Multiple weight configurations for $D$ and FNN were defined in Table~\ref{tab:weightings-D-FNN-20} to observe how the learning weights assigned to the two objective functions affect the diversity and accuracy of the generated variables.

\begin{table}[tb] 
\caption{Learning weighting definitions of $D$ and FNN}
\label{tab:weightings-D-FNN-20}%
\begin{tabular}{@{}cccccccc@{}}
\toprule
Case & 1 & 2 & 3 & 4 & 5 \\
\midrule
FNN  & 0.3    & 0.4    & 0.5    & 0.6    & 0.7  \\
$D$  & 0.7    & 0.6    & 0.5    & 0.4    & 0.3  \\   
\botrule
\end{tabular}
\end{table}

\begin{figure}[tb]
    \centering
    \includegraphics[width=0.5\textwidth]{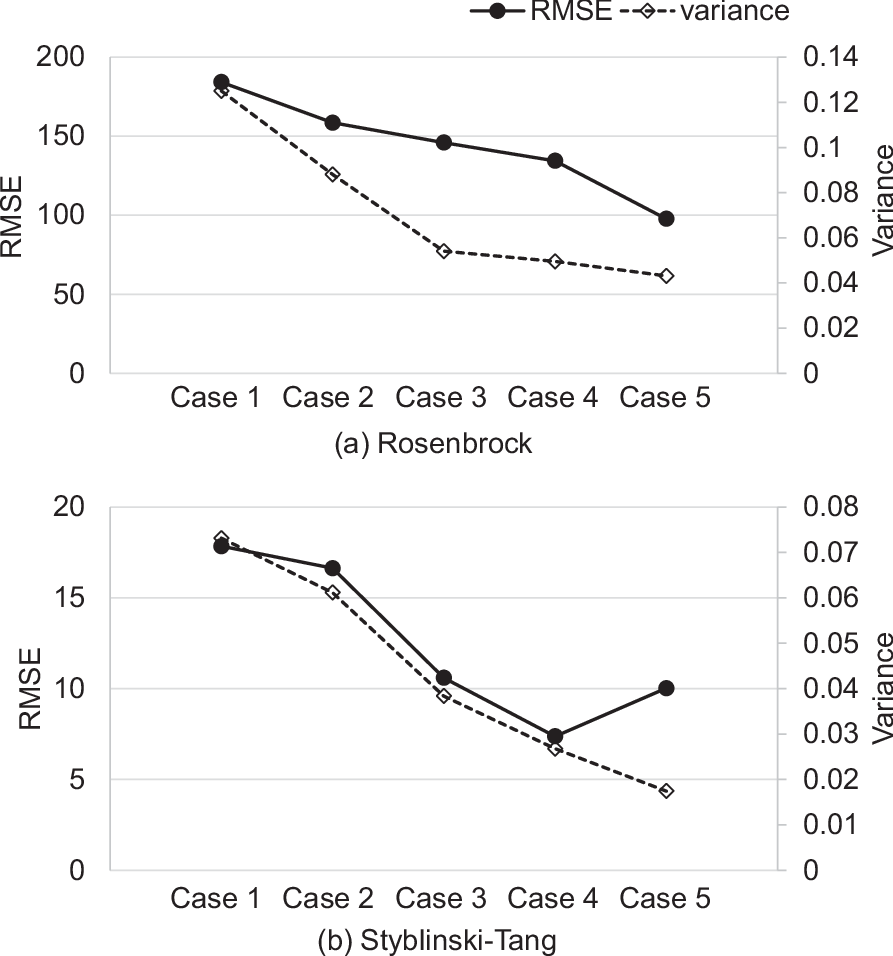}
    \caption{The optimization performance of FCGAN for the benchmark functions with different loss weights for $D$ and FNN defined in Table~\ref{tab:weightings-D-FNN-20}}
    \label{fig:FCGAN-perf-37}
\end{figure}

Fig.~\ref{fig:FCGAN-perf-37} shows the RMSE values and variance to measure the accuracy and diversity of FCGAN for each weight configuration. 
As the weight for FNN is increased, the diversity of design candidates decreases while the prediction accuracy increases, as this promotes the model to learn of the relationship between $\mathbf{x}$ and $y$. 
As the weight for $D$ is increased, the learning capability of $G$, which generates $\hat{\mathbf{x}}$, improves. 
Thus, it is possible to control the accuracy and diversity of $\hat{\mathbf{x}}$ that are typically inversely proportional by adjusting the weights of the objective functions for $D$ and FNN.
Also, one may generate design candidates with desired accuracy and diversity using unsupervised inverse design models with proper weight configurations. 
Such models may serve as a unique design optimization method that allows designers to generate diverse design candidates while achieving good accuracy.

\subsection{Deep Inverse Design with Boundary Constraints}

\begin{figure}[tb]
    \centering
    \includegraphics[width=0.5\textwidth]{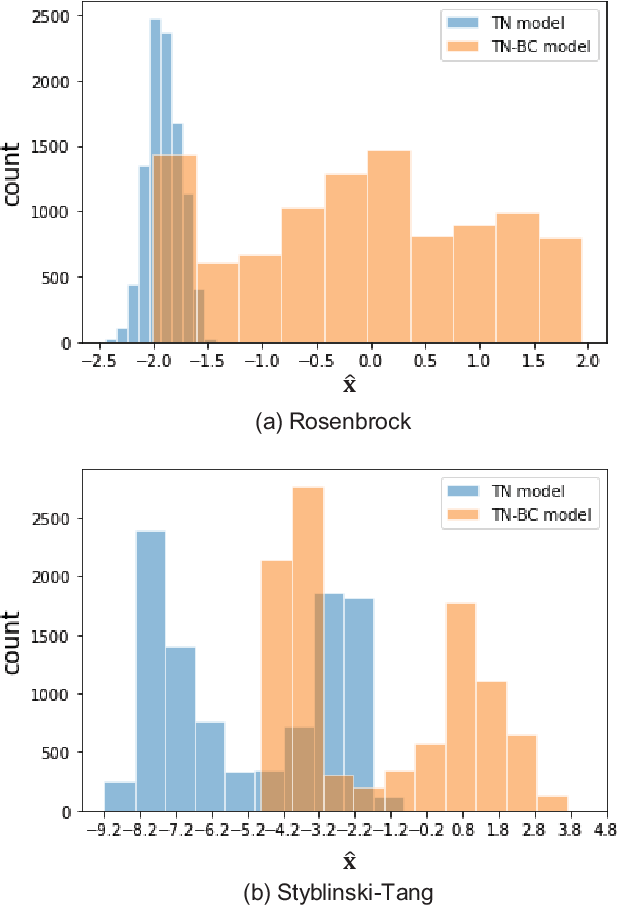}
    \caption{The histograms of $\hat{\mathbf{x}}$ from TN and TN-BC}
    \label{fig:TN-BC-hist}
\end{figure}

\begin{figure}[tb]
    \centering
    \includegraphics[width=0.5\textwidth]{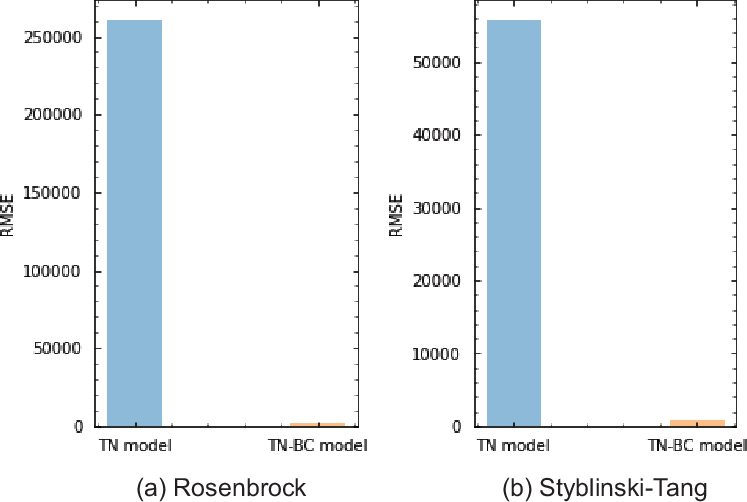}
    \caption{The RMSEs between $y$ and $y_{real}$ calculated from $\hat{\mathbf{x}}$ provided by TN and TN-BC}
    \label{fig:TN-BC-RMSE}
\end{figure}

\begin{equation}
    y_{real} = f(\hat{\mathbf{x}})
    \label{eq:y_real}
\end{equation}

The performance of TN and TN-BC was measured and compared using $y$ and $y_{real}$ derived using the ground truth benchmark function ($f$) from $\hat{\mathbf{x}}$ as shown in Eq.~\eqref{eq:y_real}.
Fig.~\ref{fig:TN-BC-hist} shows that many $\hat{\mathbf{x}}$ generated by the baseline TN was outside the boundary constraints of both benchmark functions defined in Table~\ref{tab:benchmark-data-info-2}.
The proportions of $\hat{\mathbf{x}}$ outside the ranges are 16.99\% and 51.55\% for Rosenbrock and Styblinski-Tang function, respectively.
Consequently, the gaps between $y_t$ and $y_{real}$ calculated from $\hat{\mathbf{x}}$ was much larger for the baseline TN as shown in Fig.~\ref{fig:TN-BC-RMSE}.
On the other hand, TN-BC successfully generated $\hat{\mathbf{x}}$ within the boundaries and show smaller errors.
Thus, TN-BC proves that the application of boundary constraints helps to assure good optimization performance.

\subsection{Deep Inverse Design with Various Data Compositions}

The performance of deep learning models depends on both the quality and quantity of data. 
However, it is difficult in practice to collect a large amount of good data due to time and cost constraints. 
As an alternative, it is possible to collect additional augmented data by sampling from the distribution of a small number of real observation data. 
To analyze the practicality of deep inverse design, the optimization performance was analyzed for different data size and ratio of augmentation.  
The following factors were considered:

\begin{itemize}
    \item Total data (real + augmented) sizes ($n$): 300, 500, 1,000, 4,000, 7,000, 10,000
    \item Real data ratio ($\eta$): 0\%, 25\%, 50\%, 75\%, 100\%
\end{itemize}

For each observation data ratio, the data was first collected from the observation data, and then the augmented data was sampled for the remaining ratio through non-replaceable extraction from the distribution of the observation data. 
When $\eta = 0\%$, data were extracted from a continuous uniform distribution between the minimum and maximum values of the observation data. 

\begin{figure}[tb]
    \centering
    \includegraphics[width=0.5\textwidth]{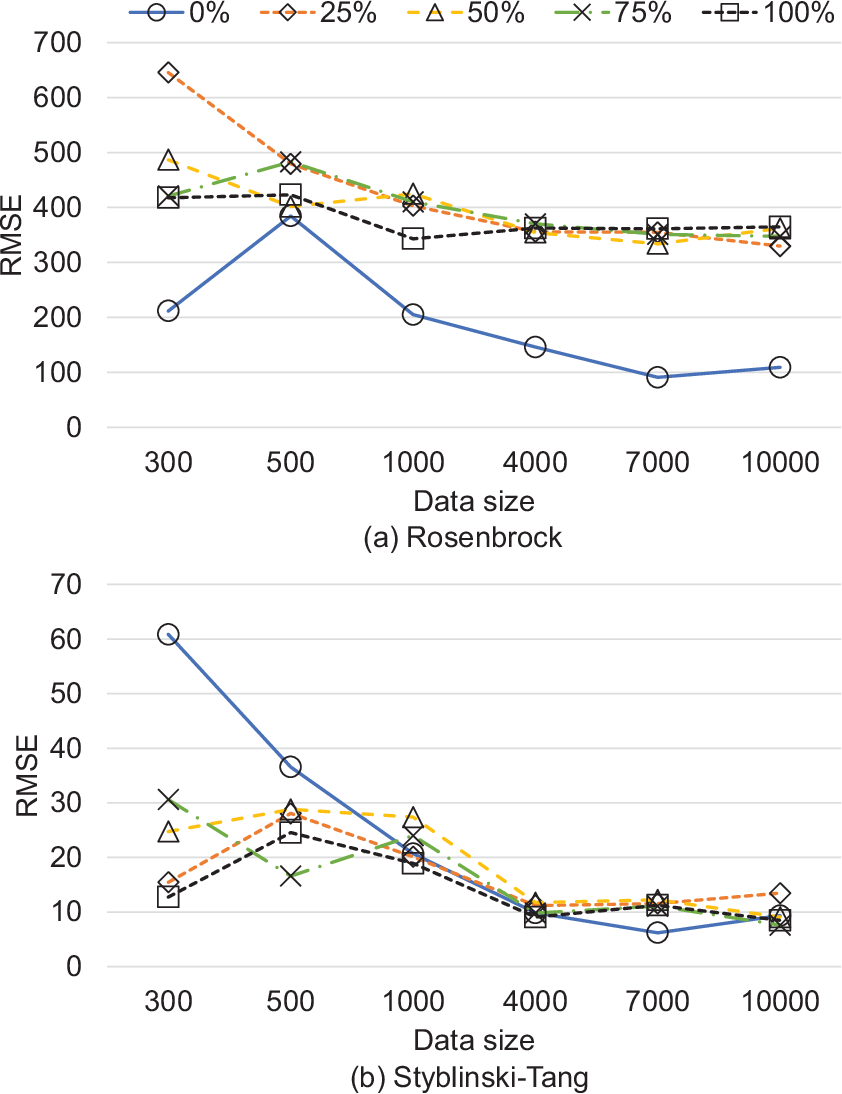}
    \caption{The optimization performance of TN with different data compositions}
    \label{fig:TN-opt-perf-data-comp-41}
\end{figure}

The TN models were trained with different data compositions defined above.
Fig.~\ref{fig:TN-opt-perf-data-comp-41} shows their performance on the test set.
The accuracy of the models improved as the data size increased, but the rate of improvement decreased once more than 1,000 data samples were used for training. 
For Rosenbrock function, the highest accuracy was achieved when trained solely on purely augmented data, when $\eta=0\%$. 
In other cases, higher stability was observed with higher ratios of real data. When $n < 500$, the performance of TN improved with higher ratios of observation data.
However,when $n > 500$, the performance was not affected by the ratio. 
In the case of Styblinski-Tang function, the worst performance was observed when trained only on pure augmented data. 
When the $n < 1,000$, a higher $\eta$ improves the performance, but the performance was not greatly affected beyond that point, similar to Rosenbrock function.

Based on this analysis, the following conclusions were drawn.
Firstly, a large amount of data is required to ensure higher stability for complex functions.
Secondly, the importance of real observation data outweighs that of augmented data. 
Lastly, once a certain number of data points are collected, the model performance becomes stable and is no longer largely affected by the data size.

\section{Conclusion}\label{sec:conclusion}

This study analyzed and compared the performance of traditional design optimization and deep inverse design methods using two benchmarks with different characteristics and case studies, assuming the surrogate model-based optimization process performed for engineering design. 
The key findings from this study can be summarized as follows.
First, deep inverse design is a promising method that may replace traditional design optimization by overcoming its limitations in terms of computational costs, while showing comparable or even higher accuracies in many cases. 
It also showed good stability regardless of dimensions or the size of a search space, showing its potential for many industrial fields with high-dimensional cases.
Second, unsupervised inverse design provides good design diversity and accuracy, and is a unique method that may balance these two as desired, that are typically inversely proportional, by adjusting their weights when training a model.
Additionally, when dealing with a high-dimensional input space $(\mathbf{x})$, it is advisable to utilize supervised learning-based inverse design rather than unsupervised learning.
Third, the use of certain activation functions may improve the optimization performance of deep inverse design models.
It was observed that when a predicted design candidate ($\hat{\mathbf{x}}$) for the target performance ($y$) satisfies the boundary constraints, the actual performance ($\hat{y}$) is closer to $y$. 
This method may significantly reduce the trial-and-error process when training deep inverse design models.
Lastly, for real-world design tasks where collecting a large amount of observation data is challenging, a specific threshold of data size may exist where the models achieve good stability.
Extracting augmented data from a distribution of small observation data for model training may also be useful.

The following future work are needed to further improve the field of deep inverse design.
First, the mathematical functions that can sample uniform data were defined as the benchmark functions to evaluate models in this study. 
However, this is typically not possible for the data used in real-world engineering often with many biases and outliers. 
Further studies on deep inverse design with more realistic benchmark functions may be necessary.
Second, the most fundamental characteristic and limitation of deep learning models is their reliance on data training. 
Recently, Physics-Informed Neural Networks (PINN) have been proposed to enable learning with simple physical information that improves the model performance with limited data \citep{daw2020physics}, which may be a promising method in terms of their applicability in the field.
Lastly, while this paper has demonstrated the innovative performance of deep inverse design models in general, further research must be conducted to examine their optimization performance in other various engineering design domains to further promote their applications.


\bibliography{sn-bibliography}

\end{document}